\documentclass[twocolum,12pt,headings=small,DIV13]{scrartcl}

\usepackage[english]{babel}
\usepackage[T1]{fontenc}
\usepackage[utf8]{inputenc}
\usepackage{mathptmx}
\usepackage{graphicx}
\usepackage{microtype}
\DeclareMathAlphabet{\mathcal}{OMS}{cmsy}{m}{n}
\usepackage[babel]{csquotes}
\usepackage{amsmath,amssymb,amsfonts,amsthm}
\usepackage{mathtools}
\usepackage{todonotes}
\usepackage[numbers,sort&compress]{natbib}

\newcommand*{\quadric}[1][Q]{\mathcal{#1}}
\newcommand*{\conic}[1][C]{\quadric[#1]}
\newcommand*{\absconic}{\conic[N]}
\newcommand*{\fregierconic}{\conic[F]}
\newcommand*{\R}{\mathbb{R}}

\newcommand{\cj}[1]{\overline{#1}}

\newtheorem{theorem}{Theorem}
\newtheorem{proposition}{Proposition}
\newtheorem{corollary}{Corollary}
\theoremstyle{definition}
\newtheorem{definition}{Definition}
\theoremstyle{remark}
\newtheorem{remark}{Remark}

\setkomafont{title}{\normalfont\bfseries}
\setkomafont{sectioning}{\normalfont\bfseries}
\setcapindent{2em}

\title{Singular Frégier Conics\\in Non-Euclidean Geometry}
\author{Hans-Peter Schröcker\\
  University of Innsbruck, Austria}

\usepackage{hyperref}
\hypersetup{
  pdfauthor={Hans-Peter Schröcker},  pdftitle={Singular Frégier Conics in Non-Euclidean Geometry},  pdfsubject={The 17th International Conference on Geometry and Graphics, Beijing, August 4-8 2016},  pdfkeywords={Frégier point, Frégier conic, Thales' theorem, hyperbolic geometry, singular conic},  hidelinks=true,
}

\begin{document}

\twocolumn[\begin{@twocolumnfalse}
  \maketitle
  \begin{abstract}
        The hypotenuses of all right triangles inscribed into a fixed conic $\conic$ with fixed right-angle vertex $p$ are incident with the Frégier point $f$ to $p$ and $\conic$. As $p$ varies on the conic, the locus of the Frégier point is, in general, a conic as well. We study conics $\conic$ whose Frégier locus is singular in Euclidean, elliptic and hyperbolic geometry. The richest variety of conics with this property is obtained in hyperbolic plane while in elliptic geometry only three families of conics have a singular Frégier locus.
   \end{abstract}
  \par\bigskip
\end{@twocolumnfalse}]

\section{Introduction}
\label{sec:introduction}

Many theorems of Euclidean elementary geometry have their counterparts in elliptic or hyperbolic geometry, possibly after a suitable re-formulation. An example of this is a version of Pythagoras' Theorem \cite{familiari-calapso69,maraner10}. By contrast, we know of no convincing non-Euclidean version for Thales' Theorem or its converse (compare \cite{weiss08}). This paper will not provide one either but we will present some interesting geometric configuration that are at least reminiscent of Thales' classical theorem.

Thales Theorem talks about right triangles with the same hypotenuse and implies that the hypotenuses of all right triangles inscribed into a circle $\conic$ contain the circle center. In particular, the hypotenuses of all right triangles with right angle vertex $p$ fixed on that circle all pass through one fixed point $f$.  This statement remains true if the circle is replaced by an arbitrary regular conic (\emph{Frégier's Theorem}). In this situation, the point $f$ is called the \emph{Frégier point} to $\conic$ and $p$.  It is not difficult to prove Frégier's theorem by means of basic projective geometry.  We present two well-known proofs: One has the benefit to clearly demonstrate the relation between Thales' and Frégier's theorems. The other employs the theory of projective transformations on conics and immediately implies the validity of Frégier's Theorem in elliptic and hyperbolic geometry.

The Frégier point $f$ depends on $\conic$ and $p$ but the locus $\fregierconic$ of all Frégier points for varying $p$ only depends on $\conic$. This locus turns out to be a conic which, by construction, shares the symmetry group with $\conic$. In Euclidean geometry, $\conic$ and $\fregierconic$ are even similar. We are interested in regular conics $\conic$ whose Frégier conic $\fregierconic$ is singular. In Euclidean geometry, these are circles and right hyperbolas. In the former case, $\fregierconic$ degenerates to the circle center, in the latter, it is the line at infinity. The same question in elliptic and hyperbolic geometry calls for a more involved answer and gives rise to a number of interesting geometric configurations, even when viewed with the eyes of a Euclidean observer.

\section{Frégier Conics in Euclidean\\geometry}
\label{sec:euclidean}

We begin by recalling a few well-known results and proofs on Frégier points and conics in the Euclidean plane. They introduce some basic concepts and set standards that later will be compared with the non-Euclidean situation.

\begin{theorem}[Frégier]
  \label{th:fregier}
  Given a regular conic $\conic$ in the Euclidean plane and a point $p \in \conic$, the hypotenuses of all right triangles inscribed into $\conic$ and with right angle at $p$ intersect in a common point~$f$.
\end{theorem}

\begin{definition}
  The point $f$ of \autoref{th:fregier} is called the \emph{Frégier point} of $\conic$ and~$p$.
\end{definition}

\begin{figure}
  \centering
  \includegraphics{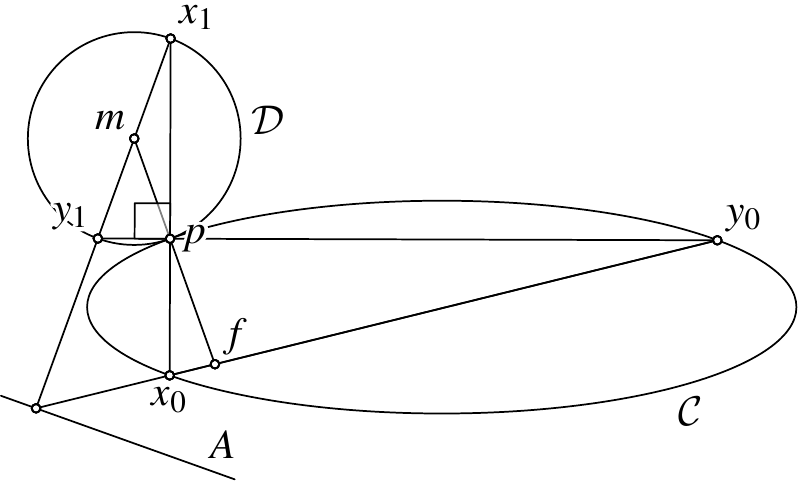}
  \caption{Proof of \autoref{th:fregier} by homology to circle}
  \label{fig:homology}
\end{figure}

We present two proofs of \autoref{th:fregier}, both having their own merits. The first proof shows how to derive Frégier's Theorem from Thales' Theorem by means of a homology to a circle (\autoref{fig:homology}).

\begin{proof}[First proof of \autoref{th:fregier}]
  Take an arbitrary circle $\conic[D]$, tangent to $\conic$ at $p$.  There exist a homology $\eta$ with center $p$ that maps $\conic[D]$ to $\conic$.  (Its axis $A$ is the Desargues axis of two triangles that correspond in $\eta$ and are inscribed into $\conic[D]$ and $\conic$, respectively.) By Thales' Theorem, the Frégier point is then $f = \eta(m)$ where $m$ is the circle center.
\end{proof}

\begin{proof}[Second proof of \autoref{th:fregier}]
  For a right triangle inscribed into $\conic$ and with right angle at $p$, denote the other vertices by $q$ and $r$.  The map $\varphi\colon \conic \to \conic$, $q \mapsto r$ (with appropriate conventions if $p$ coincides with $q$ or $r$) projects to the orthogonal involution in the line bundle around $p$.  Hence, it is an involution in $\conic$ and there exists a point, the Frégier point $f$, which is collinear with all pairs of corresponding points \cite[Theorem~8.2.8]{casas-alvero14}.
\end{proof}

This second proof appeals to a more profound knowledge of projective geometry but has the benefit of retaining its validity in non-Euclidean geometries:

\begin{corollary}
  Frégier's Theorem is true in the elliptic and hyperbolic plane.
\end{corollary}

\begin{figure}
  \centering
  \includegraphics{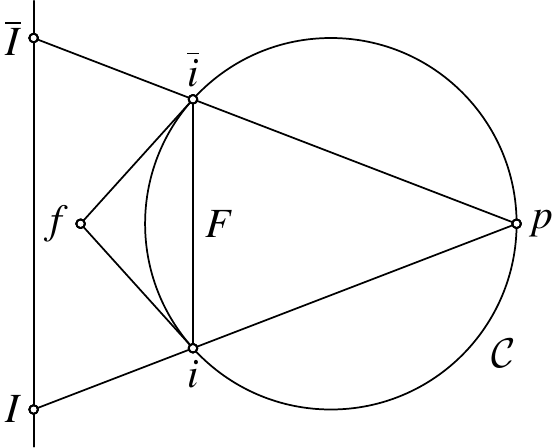}
  \includegraphics{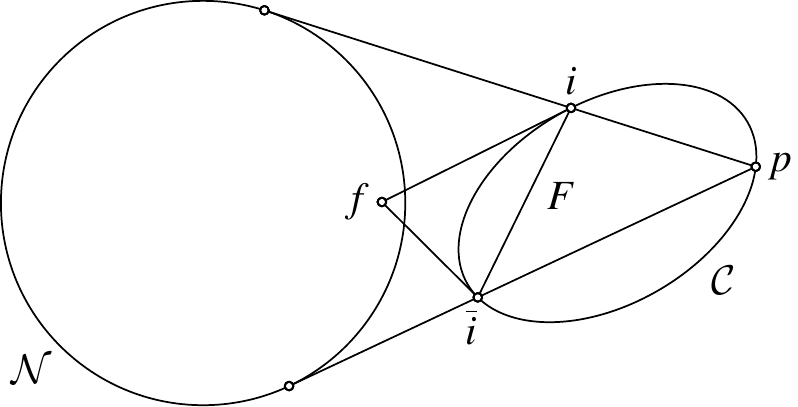}
  \caption{Construction of the Frégier point via the isotropic
    lines through~$p$ (pseudo-Euclidean and hyperbolic geometry)}
  \label{fig:involutions}
\end{figure}

The defining property of the Frégier point also allows its computation. However, the arbitrariness of the inscribed right triangle is somewhat awkward. We therefore look for alternatives. A useful observation, yet insufficient to nail down $f$, is the fact that it is located on the conic normal at $p$. As suggested in \cite[Exercise~10.12]{casas-alvero14}, we are led to consider the fix points of the involution $\varphi$ in our second proof of \autoref{th:fregier}. They are the intersection points $i$, $\cj{i}$ different from $p$ of $\conic$ and the isotropic lines through $p$. Their tangents intersect in $f$. These arguments equally apply to Euclidean and non-Euclidean geometries. \autoref{fig:involutions} displays this construction in the pseudo-Euclidean plane with absolute points $I$, $\cj{I}$ and in the hyperbolic plane with absolute conic~$\absconic$.

\begin{proposition}
  \label{prop:fregier-construction}
  The Frégier point $f$ to a conic $\conic$ and a point $p \in \conic$ in Euclidean, elliptic or hyperbolic geometry is the pole of the line $i \vee \cj{i}$ where $i$ and $\cj{i}$ are the projections of $p$ onto $\conic$ via the isotropic lines through~$p$.
\end{proposition}

\begin{remark}
  At this point, a remark on our view on hyperbolic geometry seems appropriate. In the sense of Wildberger's universal hyperbolic geometry \cite{wildberger13,wildberger10,wildberger11}, we treat points in- and outside the absolute conic $\absconic$ equally. This leads to simplified statements and computations. Arguably, the resulting theory is richer and more comprehensive. It also allows a ``real'' \autoref{fig:involutions} to illustrate \autoref{prop:fregier-construction}. Note however that universal hyperbolic geometry fails to be a model for the axiomatic geometry obtained by replacing Euclid's parallel postulate with its hyperbolic counterpart.
\end{remark}

Denoting the Frégier point to conic $\conic$ and point $p$ by $f(\conic,p)$, we call the set
\begin{equation*}
  \fregierconic \coloneqq \{f(\conic,p) \mid p \in \conic\}
\end{equation*}
the \emph{Frégier locus} of~$\conic$. In general, it is a regular conic section but exceptions may occur.

\begin{theorem}
  \label{th:fregier-conic}
  Generically, the Frégier locus of a conic $\conic$ is a conic~$\fregierconic$.
\end{theorem}

It is our aim in this paper to characterize regular conics $\conic$ whose Frégier locus $\fregierconic$ is not a regular conic. In doing so, we will derive the algebraic equation of the Frégier locus for different relative positions of $\absconic$ and $\conic$ and thus prove \autoref{th:fregier-conic}. This discussion also makes the word ``generically'' precise. If the Frégier locus is a conic, we call it the \emph{Frégier conic} to~$\conic$. Otherwise, we speak of the \emph{singular Frégier locus.}

\section{Singular Frégier Loci}
\label{sec:fregier-singular}

Now we have a closer look at conics in the Euclidean, elliptic and hyperbolic plane whose Frégier locus is singular. The Euclidean discussion is straightforward and yields known results. The non-Euclidean discussion requires the distinction between relative position of $\absconic$ and $\conic$, that is, different types of pencils of conics. We shall see that singular Frégier loci are of greater variety in a non-Euclidean setting.
 
\subsection{Frégier Conics in the Euclidean Plane}

Using homogeneous coordinates $[x_0,x_1,x_2]$, an ellipse or hyperbola in the Euclidean plane can always be described by an equation of the shape
\begin{equation*}
  \conic\colon bx_1^2 + ax_2^2 = x_0^2
\end{equation*}
with non-zero real numbers $a$, $b$ that are not both negative. The Frégier conic, computed via \autoref{prop:fregier-construction}, has equation
\begin{equation*}
  \fregierconic\colon b(a+b)^2x_1 + a(a+b)^2x_2 = ab(a-b)^2x_0.
\end{equation*}
We see that, in general, $\fregierconic$ is obtained from $\conic$ by a scaling with factor $(a-b)/(a+b)$ about the common center of $\conic$ and $\fregierconic$. This statement is not true if the conic $\conic$ is a circle ($a^2-b^2=0$) or a right hyperbola ($a^2+b^2=0$). In the former case, the Frégier conic consists of a single point, in the latter, the Frégier conics degenerates to the line at infinity. It is noteworthy that in this case the map $p \in \conic \to f \in \fregierconic$ is a double cover of the line at infinity but only the ideal points on normals of $\conic$ arise as real Frégier points. Thus, the Frégier locus is a \emph{projective line segment}.

A parabola in the Euclidean plane may be described by the equation $\conic\colon x_0x_2 = ax_1^2$ with $a \in \R \setminus \{0\}$. Its Frégier conic $\fregierconic\colon x_0x_2 = ax_1^2+\tfrac{1}{a}x_0^2$ is just a translate of $\conic$ and we can summarize:

\begin{proposition}
  If the Frégier locus in the Euclidean plane is not a regular conic then either $\conic$ is a circle and $\fregierconic$ is its center or $\conic$ is a right hyperbola and $\fregierconic$ consists of those ideal points that belong to normal directions of~$\conic$.
\end{proposition}

\subsection{Frégier Conics in the Hyperbolic Plane}

Our investigation of Frégier conics in the hyperbolic plane is based on a discussion of the relative position of the absolute conic $\absconic$ and the conic $\conic$ in the \emph{complex} projective plane, that is, pencils of conics in that plane. This is justified because the Frégier locus to a real conic is always real, even if some elements in the construction of \autoref{fig:involutions} appear as conjugate imaginary pairs. The line $F$ and its pole $f$ are always real.

\begin{figure*}
  \centering
  \includegraphics{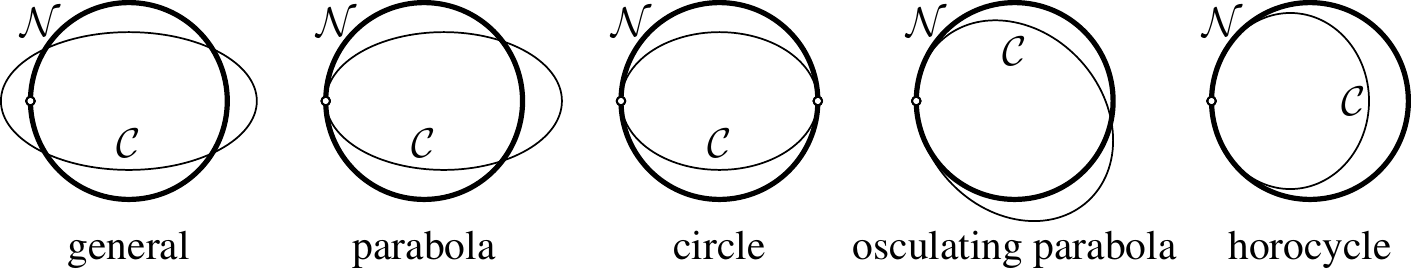}
  \caption{Pencils of conics, conics in the hyperbolic plane.}
  \label{fig:pencils}
\end{figure*}

The conics of a pencil share four different points (``base points''), some of which may coincide and thus result in tangency or contact of higher order \cite[Section~9.6]{casas-alvero14}. Depending on the number of coinciding points, one can distinguish five different cases (\autoref{fig:pencils}):
\begin{enumerate}
\item General pencils with four different base points.
\item Simple contact pencils with a single pair of coinciding base points. In analogy to the Euclidean situation, we call the corresponding conics \emph{parabolas.}
\item Bitangent pencils where precisely two pairs of base points coincide. We call the conics in this case \emph{circles.}
\item Double contact or osculating pencils where three of the four base points coincide. Here, we speak of \emph{osculating parabolas.}
\item Triple contact or hyperosculating pencils where all four base points coincide. The corresponding conics are called \emph{horocycles.}
\end{enumerate}

\subsubsection{General Conics}

The equation of a general conic $\conic$ may be written as
\begin{equation*}
  \conic\colon bx_1^2 + ax_2^2 = x_0^2,
  \quad
  a, b \in \R \setminus\{0\}
\end{equation*}
with non-zero real numbers $a$, $b$ that are not both negative and not
both equal to $1$. This case also comprises circles for $a = 1$,
$b \neq 1$ or $b = 1$, $a \neq 1$.

\begin{figure}
  \centering
  \includegraphics{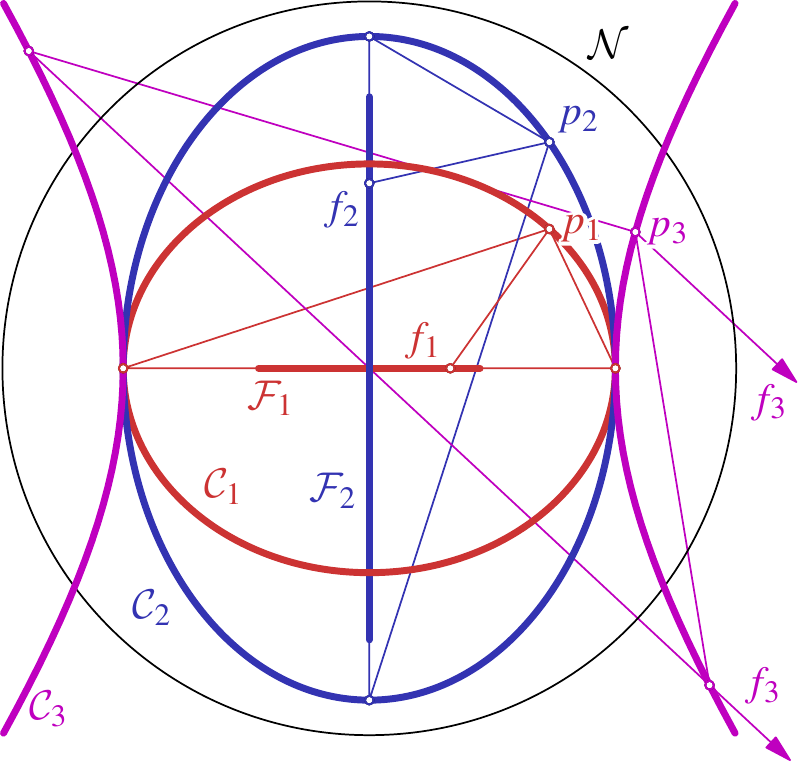}
  \caption{Three conics with singular Frégier locus ($\fregierconic_3$
    is part of the line at infinity)}
  \label{fig:fregierConicEllHyp}
\end{figure}

Its Frégier conic has the equation
\begin{multline*}
  \fregierconic\colon \frac{(a^2b^2-a^2-b^2)^2}{(a^2b^2+a^2-b^2)^2a^2}x_1^2\\
  + \frac{(a^2b^2-a^2-b^2)^2}{(a^2b^2-a^2+b^2)^2b^2}x_2^2 = x_0^2.
\end{multline*}
It is singular if and only if
\begin{equation*}
  b^2 = \frac{a^2}{a^2+1},\
  b^2 = \frac{-a^2}{a^2-1},\text{ or }
  b^2 = \frac{a^2}{a^2-1}.
\end{equation*}
The Frégier locus is, in that order, a projective line segment on the
first, the second or the third axis of the underlying projective
coordinate frame (\autoref{fig:fregierConicEllHyp}).

\subsubsection{Parabolas}

\begin{figure}
  \centering
  \includegraphics{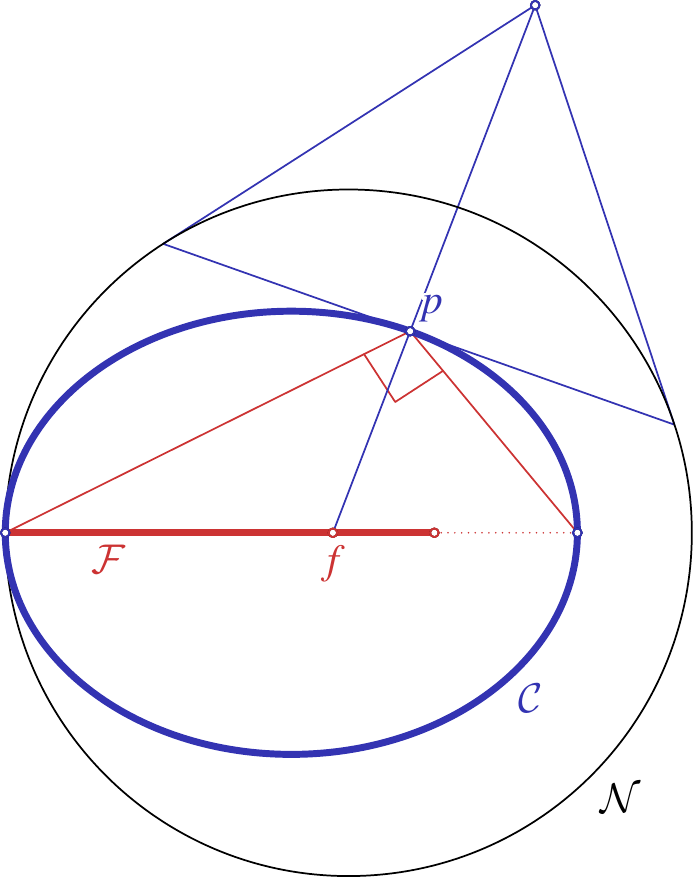}
  \caption{Parabola with singular Frégier locus.}
  \label{fig:fregierConicParabola}
\end{figure}

A parabola that is tangent to $\absconic$ at $[1,-1,0]$ admits an
equation of the shape
\begin{multline*}
  \conic\colon \lambda(x_1^2+x_2^2-x_0^2)\\+ (x_0+x_1)(\mu(x_0+x_1)+x_1) = 0.
\end{multline*}
The parameters $\lambda$ and $\mu$ range in $\R$ but $\lambda = 0$ and
$\lambda = -\tfrac{1}{2}$ are prohibited in order to ensure regular
conics. The limiting case $\mu \to \infty$ yields a horocycle and will
be treated later. The Frégier conic has equation
$\fregierconic\colon (x_1^2+x_2^2-x_0^2)\lambda^4 + (x_0+x_1)\Lambda$
where
\begin{multline*}
  \Lambda = (x_0+x_1)(5\mu\lambda^3+12\mu\lambda^2+9\mu\lambda+2\mu+1)\\+
  \lambda(4(\lambda+1)x_0 + (5\lambda^2+8\lambda+5)x_1).
\end{multline*}
It is singular precisely for $\lambda = 0$, $\lambda = -\tfrac{1}{2}$,
or $\lambda = -1$. Because only the last value is admissible, we obtain a
one-parametric family of hyperbolic parabolas with singular Frégier
conic:
\begin{gather*}
  \conic\colon x_0^2-x_1^2-x_2^2 + (x_0+x_1)(\mu(x_0+x_1)+x_1) = 0
\end{gather*}
and $\mu$ ranges in~$\R$ (\autoref{fig:fregierConicParabola}).

\subsubsection{Osculating Parabolas}

Writing the equation of an osculating parabola as
\begin{equation*}
  \conic\colon
  \lambda(x_1^2+x_2^2-x_0^2) + (x_0+x_1)x_2 = 0,
\end{equation*}
with $\lambda \in \R\setminus\{0\}$, the Frégier conic becomes
\begin{multline*}
  \fregierconic\colon 2\lambda(x_1^2+x_2^2-x_0^2)\\
  + (x_0+x_1)(10x_2\lambda + 8(x_0+x_1)) = 0.
\end{multline*}
It is never singular.

\subsubsection{Circles}

\begin{figure}
  \centering
  \includegraphics{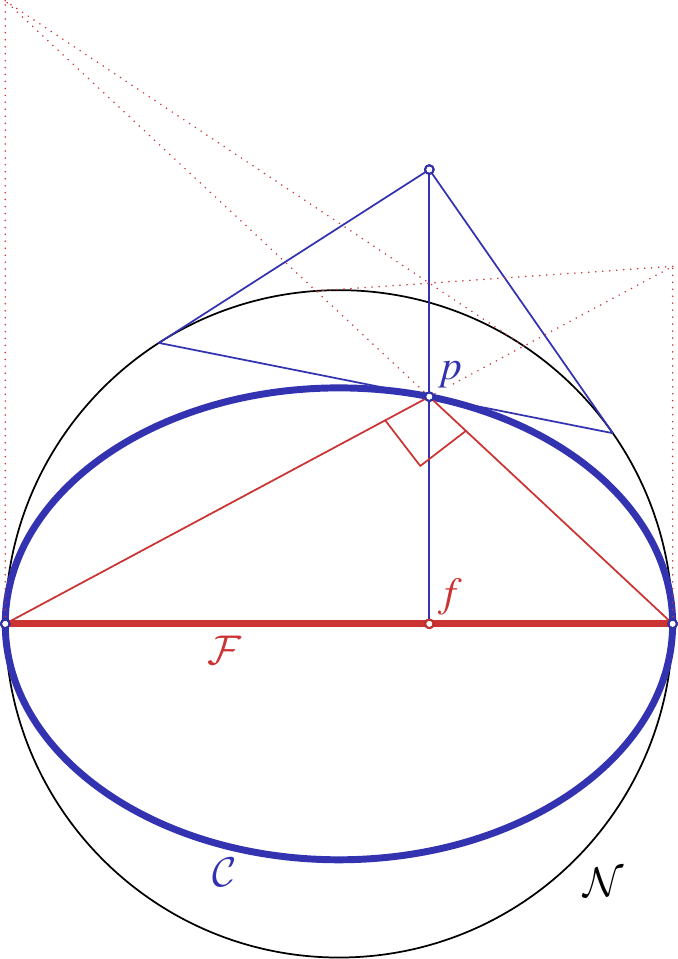}
  \caption{Singular Frégier locus to a hyperbolic circle}
  \label{fig:fregierConicCircle2-static}
\end{figure}

Circles in hyperbolic geometry are characterized by having double
contact with $\absconic$. The points of contact may be both real or
both conjugate complex. It is possibly to discuss these two cases at once
but it is probably easier to consider them separately. If the points of
tangency are real, we may write the circle equation as
\begin{gather*}
  \conic: \lambda (x_1^2+x_2^2-x_0^2) + (x_1-x_0)(x_1+x_0) = 0,\\
  \lambda \in \R \setminus \{0,1\}.
\end{gather*}
This is the parabola case for $\mu = -\frac{1}{2}$. The Frégier conic is
\begin{multline}
  \label{eq:1}
  \fregierconic: \lambda^3(x_1^2+x_2^2-x_0^2)\\+ (5\lambda^2+8\lambda+4)(x_1-x_0)(x_1+x_0) = 0.
\end{multline}
Unless $\lambda = -2$ (see below), it is again a circle.  Obviously, it
shares the symmetries of $\conic$. Because the signs of
$\lambda+1$ and $\lambda^3/(5\lambda^2+8\lambda+4)^{-1}+1$ agree, it
lies in the interior of $\absconic$ if and only if $\conic$ does. The
conic \eqref{eq:1} is singular if
\begin{equation*}
  \lambda^3(\lambda+2)^4(\lambda+1)^2 = 0.
\end{equation*}
Because $\lambda \notin \{0,1\}$, the Frégier locus is singular precisely for $\lambda = 2$. It degenerates to the line with equation $x_2 = 0$, that is, the span of the two points of tangency. More precisely, only interior point of $\absconic$ occur as real Frégier points (\autoref{fig:fregierConicCircle2-static}). The lines connecting any conic point $p \in \conic$ with the points of tangency are perpendicular.

\begin{remark}
  The situation of
  \autoref{fig:fregierConicCircle2-static} may also be described by
  saying that \enquote{Thales Theorem in the hyperbolic plane
    holds true for infinite line segments.}
\end{remark}

\begin{remark}
  \autoref{fig:fregierConicCircle2-static} is also remarkable from a
  Euclidean viewpoint. The ellipse $\conic$ with semi-axis ratio
  $1:1/\sqrt{2}$ is inscribed into the Thales circle $\absconic$ over
  major axis in such a way that for any point $p \in \conic$, the
  pole of the ellipse tangent in $p$ with respect to $\absconic$, the
  projection of $p$ onto the major axis and $p$ itself are collinear.
\end{remark}

\begin{figure}
  \centering
  \includegraphics{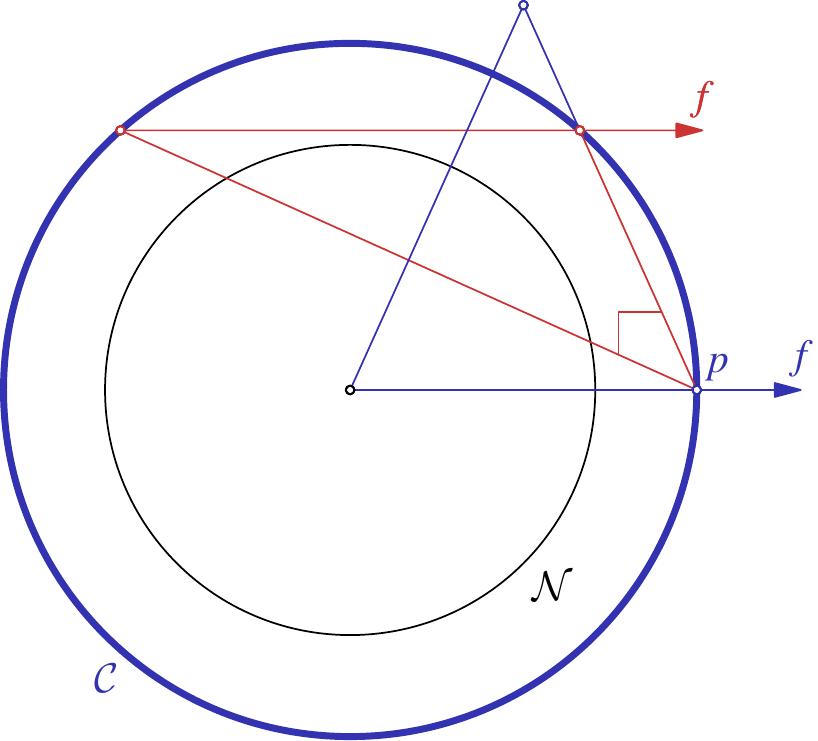}
  \caption{Singular Frégier locus to a hyperbolic circle}
  \label{fig:fregierConicCircle1-static}
\end{figure}

If both points of tangency are conjugate complex, the circle equation becomes
\begin{gather*}
  \conic: \lambda (x_1^2+x_2^2-x_0^2) + x_1^2+x_2^2 = 0,\\
  \lambda \in \R \setminus \{0,-1\}
\end{gather*}
and the Frégier conic is
\begin{multline*}
  \fregierconic: \lambda^3(x_1^2+x_2^2-x_0^2)\\
  - (5\lambda^2+8\lambda+4)(x_1^2+x_2^2) = 0.
\end{multline*}
Unless $\lambda = -2$, symmetry is shared between $\conic$ and $\fregierconic$ and one lies in the interior of $\absconic$ if and only if the other does. The Frégier conic is singular precisely for $\lambda = -2$ (\autoref{fig:fregierConicCircle1-static}) and it is the span of the two points of tangency. Viewing \autoref{fig:fregierConicCircle1-static} with the eyes of a Euclidean observer, the Frégier locus is a projective line segment on the line at infinity.

\subsubsection{Horocycles}

If the intersection points of $\conic$ and $\absconic$ all
coincide, $\conic$ is called a \emph{horocycle.} Its equation may be
written as
\begin{equation*}
  \conic: \lambda (x_1^2+x_2^2-x_0^2) + (x_2-x_0)^2 = 0,
  \quad \lambda \in \R \setminus \{0\}
\end{equation*}
whence the the Frégier conic becomes
\begin{equation*}
  \fregierconic\colon \lambda(x_1^2+x_2^2-x_0^2) + 5(x_2-x_0)^2 = 0.
\end{equation*}
We see that $\conic$ and $\fregierconic$ are related by the map
$\lambda \mapsto \frac{1}{5}\lambda$. No singularities arise.

\subsection{Frégier Conics in the Elliptic Plane}

The situation in elliptic geometry is algebraically equivalent to
universal hyperbolic geometry. Via a complex projective
transformation, all results of the latter also hold in the projective
extension of the former. In order to make statements on Frégier conics
in the elliptic plane over the real numbers, we only need to discuss
the reality of the involved geometric objects. The only relevant case
is that of a general conic and its specialization to a circle:
\begin{equation*}
  \conic\colon bx_1^2 + ax_2^2 = x_0^2,\quad
  a,b\in\R \setminus\{0\}.
\end{equation*}
Keeping things short, we only mention the final result. In the
elliptic plane, the Frégier conic to $\conic$ is singular if and only
if $a$ and $b$ satisfy one of
\begin{equation*}
  b^2 = \frac{a^2}{a^2+1},\
  b^2 = \frac{-a^2}{a^2-1},\text{ or }
  b^2 = -\frac{a^2}{a^2+1}.
\end{equation*}
The last case is never real and none of these families contains real circles.

\section{Conclusion}

We recalled some well-known facts about Frégier conics in the
Euclidean plane and transferred them to elliptic and hyperbolic
geometry. It turned out that the situation between Euclidean and
non-Euclidean geometry is different when it comes to singular Frégier
loci. In particular, we saw that the Frégier locus to a regular conic
$\conic$ can be singular if $\conic$ is an ellipse, a hyperbola, a
parabola, or a circle but not if it is an osculating parabola or a
horocycle. Another notable difference is that a singular Frégier locus
in non-Euclidean geometry is always a projective line segment while it
may be a single point in Euclidean geometry.

 \bibliographystyle{plainnat}

\begin{thebibliography}{7}
\providecommand{\natexlab}[1]{#1}
\providecommand{\url}[1]{\texttt{#1}}
\expandafter\ifx\csname urlstyle\endcsname\relax
  \providecommand{\doi}[1]{doi: #1}\else
  \providecommand{\doi}{doi: \begingroup \urlstyle{rm}\Url}\fi

\bibitem[Casas-Alvero(2014)]{casas-alvero14}
Eduardo Casas-Alvero.
\newblock \emph{Analytic Projective Geometry}.
\newblock European Mathematical Society, 2014.
\newblock ISBN 978-3-03719-138-5.

\bibitem[Familiari-Calapso(1969)]{familiari-calapso69}
Maria~Teresa Familiari-Calapso.
\newblock Sur une classe de triangles et sur le théorème de pythagore en
  géométrie hyperbolique.
\newblock \emph{C. R. Acad. Sci. Paris Sér. A--B}, 263:\penalty0 A668--A670,
  1969.

\bibitem[Maraner(2010)]{maraner10}
Paolo Maraner.
\newblock A spherical {Pythagorean} theorem.
\newblock \emph{Math. Intell.}, 32\penalty0 (3):\penalty0 46--50, 2010.

\bibitem[Weiß and Gruber(2008)]{weiss08}
Gunter Weiß and Franz Gruber.
\newblock {Den Satz von Thales verallgemeinern -- aber wie?}
\newblock \emph{KoG}, 12:\penalty0 7--18, 2008.

\bibitem[Wildberger(2010)]{wildberger10}
Norman~J. Wildberger.
\newblock Universal hyperbolic geometry {II:} a pictorial overview.
\newblock \emph{KoG}, 14\penalty0 (14):\penalty0 3--24, 2010.

\bibitem[Wildberger(2011)]{wildberger11}
Norman~J. Wildberger.
\newblock Universal hyperbolic geometry {III:} first steps in projective
  triangle geometry.
\newblock \emph{KoG}, 15\penalty0 (15):\penalty0 25--49, 2011.

\bibitem[Wildberger(2013)]{wildberger13}
Norman~J. Wildberger.
\newblock Universal hyperbolic geometry {I:} trigonometry.
\newblock \emph{Geometriae Dedicata}, 163\penalty0 (1):\penalty0 215--274,
  2013.

\end{thebibliography}

\end{document}